\documentclass[12pt]{amsart}

\usepackage{color}
\usepackage{amssymb}

\ifx\pdfoutput\undefined
  \usepackage{graphicx} 
  \DeclareGraphicsExtensions{.pstex,.eps}
\else
  \ifx\pdfoutput\relax
     \usepackage{graphicx} 
     \DeclareGraphicsExtensions{.pstex,.eps}
  \else
     \usepackage{amsfonts}
     \usepackage[pdftex]{graphicx}  \DeclareGraphicsExtensions{.pdf,.mps}
  \fi
\fi

\newcommand{\RR}{{\mathbb R}}

\renewcommand{\l}{\lambda}

\renewcommand{\Im}{\mathop{\rm Im}\nolimits}

\theoremstyle{plain}

\theoremstyle{definition}

\numberwithin{equation}{section}

\def\squarebox#1{\hbox to #1{\hfill\vbox to #1{\vfill}}} 
 
\newcommand{\sech}{\textnormal{sech}}

\newcommand{\nlso}{\textnormal{NLS}_0}
\newcommand{\nlsq}{\textnormal{NLS}_q}
\usepackage{amsxtra}

\title
[Soliton splitting by external delta potentials]
{Soliton splitting by external delta potentials}

\author[J. Holmer]
{Justin Holmer}
\author[J. Marzuola]
{Jeremy Marzuola}
\author[M. Zworski]
{Maciej Zworski}
\address{Mathematics Department, University of California \\
Evans Hall, Berkeley, CA 94720, USA}

\begin{document}

\begin{abstract}
We show that a soliton scattered by an external delta potential 
splits into two solitons and a radiation term. Theoretical 
analysis gives the amplitudes and phases of the reflected and
transmitted solitons with errors going to zero as the velocity 
of the incoming soliton tends to infinity. Numerical analysis
shows that this asymptotic relation is valid for all but very 
slow solitons. 

We also show that the total transmitted mass, that is 
the square of the $L^2$ norm of the solution restricted on the transmitted side of the delta potential  
is in good agreement with the quantum transmission 
rate of the delta potential. 

\end{abstract}

\maketitle

\section{Introduction and statement of results}   
\label{in}

A bright soliton is a travelling wave solution,
\begin{equation}
\label{eq:trwv}
u ( x , t ) =  A \sech ( A ( x - v t ) ) 
\exp (i \varphi +  i v x + i( A^2 - v^2) t/ 2 ) \,, \ \ A > 0 \,, \ \ v \in \RR \,, 
\end{equation}
of the nonlinear Schr\"odinger equation (NLS):
\begin{equation}
\label{eq:nls1}
i\partial_t u + \tfrac{1}{2}\partial_x^2 u +u|u|^2 = 0
\end{equation}
Its remarkable feature is total coherence -- see for instance \cite{KiMa} for a review of theoretical and experimental situations in which bright solitons arise. 

\begin{figure}
\scalebox{0.5}{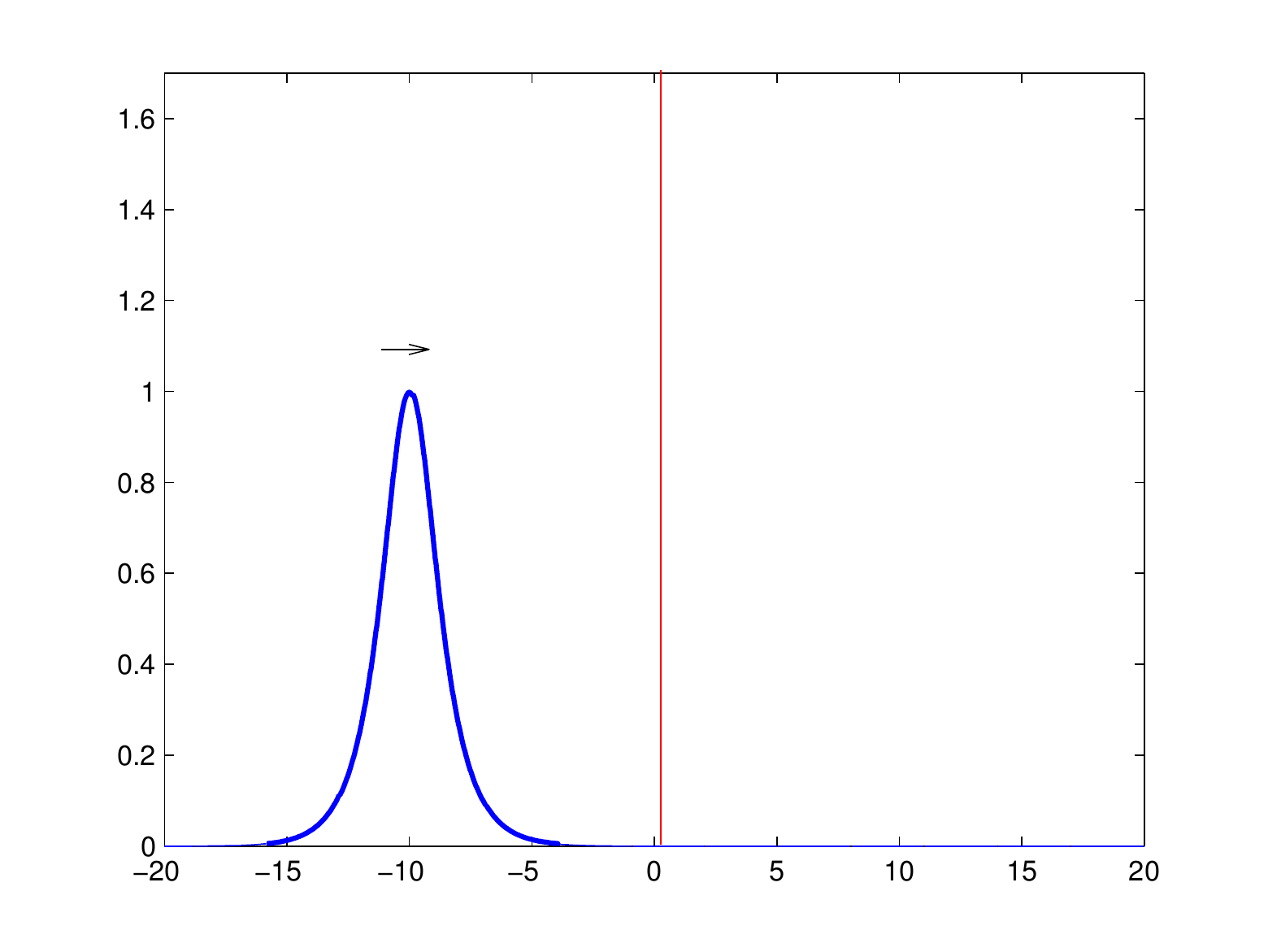} \hfill \scalebox{0.5}{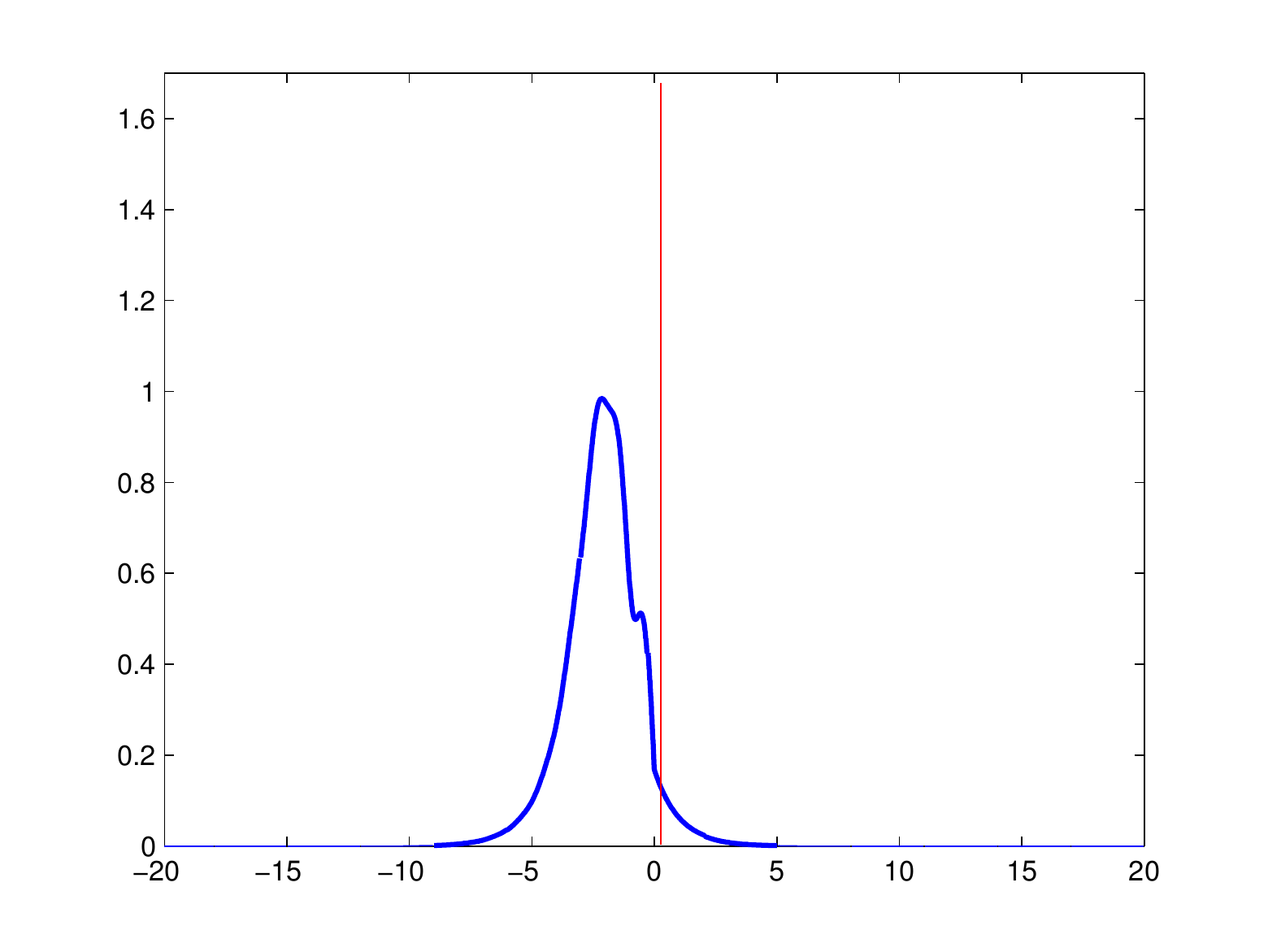}
\scalebox{0.5}{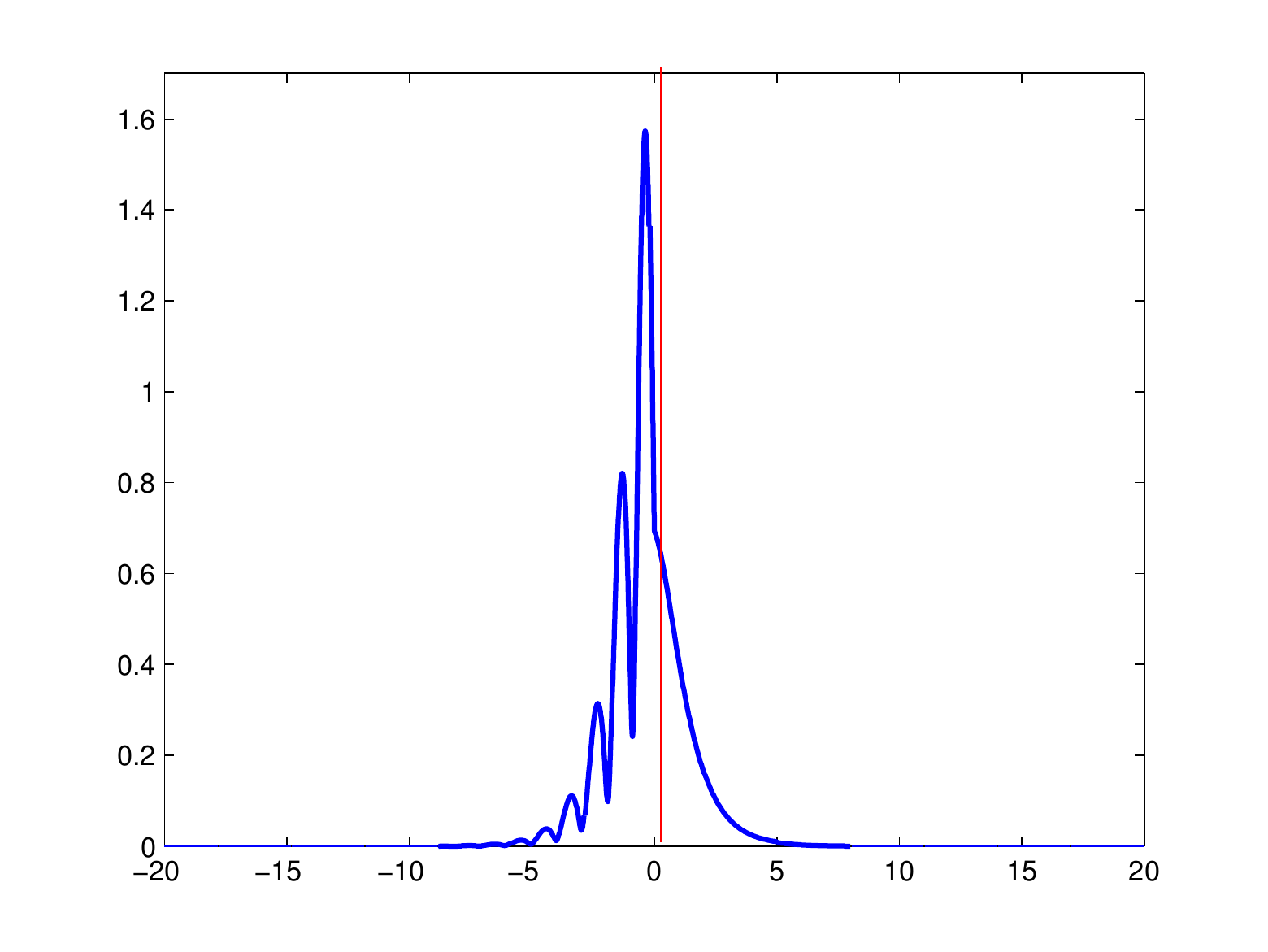} \hfill \scalebox{0.5}{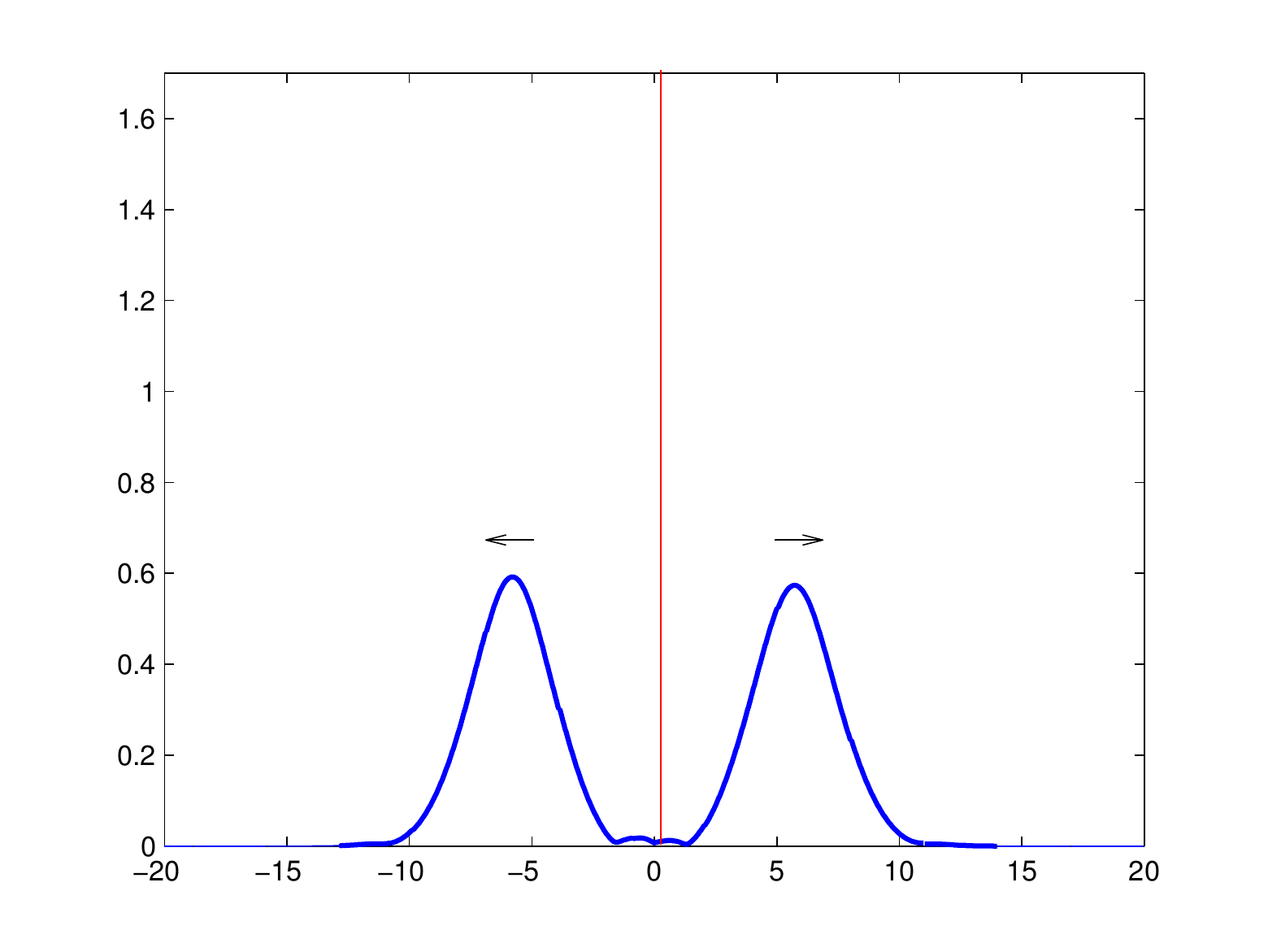}

\caption{\label{f:solv3} Numerical simulation of the case $q=v=3$, $x_0=-10$,
at times $t=0.0, 2.7, 3.3, 4.0$.  Each frame is a plot of amplitude $|u|$ versus $x$.}
\end{figure}

Suppose now that we consider a perturbed NLS, that is, the Gross-Pitaevskii equation,
by adding an external  potential:
\begin{equation}
\label{eq:nls}
\left\{
\begin{aligned}
&i\partial_t u + \tfrac{1}{2}\partial_x^2 u -q\delta_0(x)u +u|u|^2 = 0\\
&u(x,0) = u_0(x)
\end{aligned}
\right.
\end{equation}
If as 
 initial data we take a soliton approaching the impurity 
from the left:
\begin{equation}
\label{eq:init}
 u_0 ( x ) = e^{ i v x } \sech ( x - x_0 ) \,, 
\ \  x_0 \ll 0 \,, 
\end{equation}
then until time $ t_0 \approx x_0 / v $, the propagation will still be approximately given 
by \eqref{eq:trwv}. Here we put $ A = 1 $ and $ \varphi = 0$. Scaling properties of the 
delta function show that this allows general soliton initial conditions. Thus the velocity, $ v $, 
and the coupling constant, $ q $, are the only parameters of the problem.

For $ t >  x_0/v $ the effects of the delta potential are 
dramatically visible and as we show in this paper they can be understood using 
the transmission and reflection coefficients of the delta potential from 
standard scattering theory.

For the soliton scattering the natural definition of the 
transmission rate is given by 
\begin{equation}
\label{eq:tqs} 
 T_q^{\rm{s}} ( v ) = \frac{1} 2  \lim_{ t\rightarrow \infty } 
\int_0^\infty  | u ( t , x ) |^2 dx \,, \ \ \int_{\RR } | u ( t,x ) |^2 dx = 2 \,,
 \end{equation}
where on the right we recalled the conservation of the $ L^2 $ norm.
The reflection coefficient is 
\begin{equation}
\label{eq:rqs} 
R_q^{\, \rm{s}} ( v )  =  \frac{1} 2  \lim_{ t\rightarrow \infty } 
\int_{-\infty}^{\, 0 }  | u ( t , x ) |^2 dx 
\,,\end{equation}
and $ T_q^{\rm{s}} ( v ) +  R_q^{\rm{s}} ( v )  = 1 $.

\begin{figure}
\scalebox{0.7}{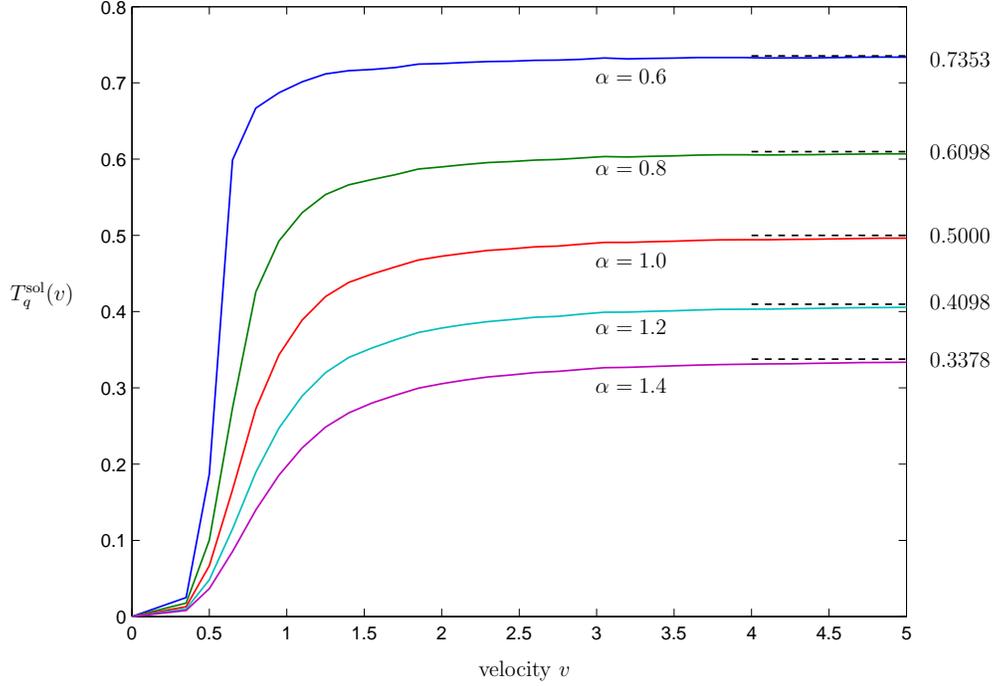}
\caption{The figure illustrates the convergence, as $v\to \infty$, of 
$T_q^{\rm{s}}(v)$ to the expected asymptotic value $1/(1+\alpha^2)$ for $\alpha=+0.6, +0.8, \ldots, +1.4$ (so $q>0$). It shows that the large velocity asymptotic behaviour in fact takes hold by velocity $v\sim 3$. }
\label{F:trans1}
\end{figure}

The following result is obtained by a numerical analysis of the problem:
\begin{equation}
\label{eq:num1}
  T^{\rm{s}}_q ( v ) = \frac{ v^2}{ v^2 + q^2 } + {\mathcal O} \left( \frac 1 {v^2   } 
\right)  \,. 
\end{equation}
In \S \ref{ros} we explain how a weaker rigorous result is obtained in 
\cite[Theorem 1]{HMZ}. Fig.\ref{F:trans1} shows the numerical agreement of 
$ T_q^{\rm{s}} ( v ) $ as a function of $ \alpha =  q / v $.

The leading term on the right hand side of \eqref{eq:num1} has the following
natural interpretation in elementary scattering theory, see for instance
\cite{LL}. Since we need it 
below to formulate the result about soliton splitting
\eqref{eq:th2} we review the basic concepts. Thus let
\[ H_q = - \frac{1}2 \frac{d^2}{dx^2} + q \; \delta_0 ( x ) \,,\]
and consider a general solution to $ ( H_q - \lambda^2 /2 ) u = 0 $,
\[  
u ( x) = A_\pm e^{ - i \lambda x } + B_\pm e^{ i \lambda x } \,, \ \pm x > 0 
\,.\]
The matrix 
\[ S ( \lambda ) \, : \,  \begin{bmatrix}
A_+ \\ B_- \end{bmatrix}  \longmapsto 
 \begin{bmatrix} A_- \\ B _+ \end{bmatrix} \,,\]
is called the {\em scattering matrix} and in our simple case it can 
be easily computed:
\[ S ( \lambda ) = \begin{bmatrix} t_q ( \lambda ) & r_q ( \lambda ) 
\\ r_q ( \lambda ) & t_q ( \lambda ) \end{bmatrix} \,,\]
where $ t_q $ and $ r_q $ are the transmission and reflection coefficients:
\begin{equation}
\label{eq:tr}
t_q ( \lambda ) = \frac{ i \lambda } { i \lambda - q } \,, \ \ 
r_q ( \lambda ) = \frac{ q} {i \lambda - q } \,.
\end{equation}
They satisfy two equations, one standard (unitarity) 
and one due to the special structure of the potential:
\begin{equation}
\label{eq:trpr} | t_q ( \lambda ) |^2 + | r_q ( \l ) |^2 = 1 \,, \ \ 
t_q ( \lambda ) = 1 + r_q ( \lambda ) \,.\end{equation}

The quantum transmission rate at velocity $ v $ is given 
by the square of the absolute value of the transmission coefficient 
\eqref{eq:tr},
\begin{equation}
\label{eq:tqv}
  T_q ( v ) = | t_q ( v ) |^2 = \frac{ v^2}{ v^2 + q^2 } \,.
\end{equation}
We recall (see \cite[(2.21)]{HMZ}) that if $ \psi $ is a smooth
function which is zero outside, say, $ [-a,-b] $, $ a > b > 0 $, 
then 
\[ \int_0^\infty | \exp ( - i t H_q ) \psi ( x ) |^2 dx = 
T_q  ( v ) + {\mathcal O} \left( \frac 1 {v^2} \right) \,, \]
just as in the nonlinear soliton experiment \eqref{eq:num1}.

Hence \eqref{eq:num1} shows that in scattering of fast solitons the transmission rate is well approximated by the quantum transmission rate of the delta potential -- see \S \ref{nr} for more on that and the comparison with the linear case.

Our second result shows that the scattered solution is given by a sum of a reflected and a transmitted soliton, and of a time decaying (radiating) term. In other words, the delta potential splits the incoming soliton into two waves which become single {\em solitons}.  In previous works in the physics literature (see for instance \cite{CM}) the resulting waves were only described as ``soliton-like''.  More precisely, for $  t \gg |x_0|/v $  we have 
\begin{gather}
\label{eq:th2}
\begin{gathered}
u(x,t)  = u_T ( x , t) + u_R ( x , t ) 
 + \mathcal{O}_{L_x^\infty}\left(\left(t-{|x_0|}/{v}\right)^{-1/2}\right) 
+ {\mathcal O}_{L_x^2}( v^{-2} )  \,,  \\
 u_T ( x , t )  =  e^{ i \varphi_T} 
e^{ixv + i  ( A_T^2 - v^2 )t /2 }  A_T \, \sech(A_T(x-x_0-tv)) \,, \\
u_R ( x, t)   =  e^{i \varphi_R} 
 e^{-ixv + i  ( A_R^2 - v^2 ) t/2 } A_R \, \sech(A_R(x+x_0+tv)) \,, 
\end{gathered}
\end{gather}
where 
\[ A_T = \left\{ \begin{array}{ll} 
2|t_q(v)|-1 \,, & |t_q ( v ) | \geq 1/2 \\
\ \ \ 0\,, & | t_q ( v ) | \leq 1/2 \,, 
\end{array} \right.  \qquad 
 A_R =  \left\{ \begin{array}{ll} 
2|r_q(v)|-1\,, & |r_q ( v ) | \geq 1/2 \\
\ \ \ 0\,,  & | r_q ( v ) | \leq 1/2 \,,
\end{array} \right.  \] 

and
\begin{gather*}
\varphi_T = \arg t_q ( v ) + \varphi_0(|t_q(v)|) + (1-A_T^2)|x_0|/2v \,, 
\\
\varphi_R = \arg r_q ( v ) + \varphi_0(|r_q(v)|) + (1-A_R^2)|x_0|/2v \,, 
\end{gather*}
$$\varphi_0(\omega) =  \int_0^\infty \log\left( 1 + \frac{\sin^2\pi \omega}{\cosh^2\pi \zeta} \right) \frac{\zeta}{\zeta^2+(2\omega-1)^2} \, d\zeta \,. $$
Here $ t_q ( v ) $ and $ r_q ( v ) $ are the transmission and reflection coefficients of the delta-potential (see \eqref{eq:tr}). 

\begin{figure}
\begin{center}
\scalebox{0.7}{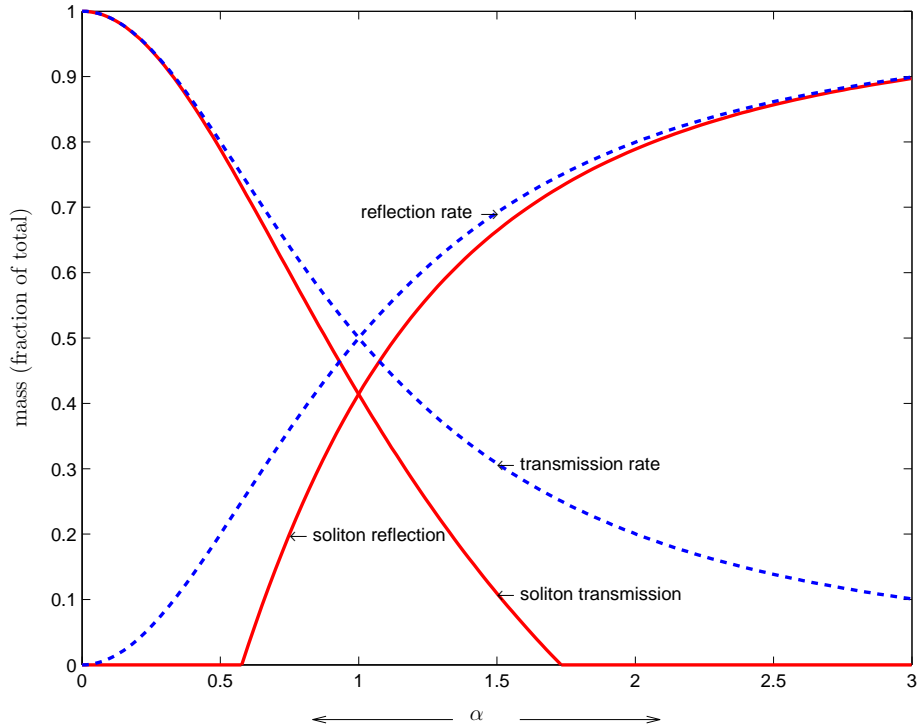}
\end{center}
\caption{A comparison of the linear and nonlinear (soliton) scattering 
rates as functions of $\alpha = q/v $.}
\label{f:new}
\end{figure}

The result is illustrated in Fig.\ref{f:new}. We can consider $ A_R ( q/v ) $ and $ A_T ( q /v ) $ as nonlinear replacements of $ R_q ( v ) $ and $ T_q ( v ) $, respectively.  Clearly $ A_T + A_R \neq 1$ except in the asymptotic limits $ q / v \rightarrow 0 \,, \infty $.  Thus if we consider soliton scattering ``particle-like'' it is nonelastic. 

In Fig.\ref{f:new} we also see the thresholds for the formation of 
reflected and transmitted solitons: 
\begin{equation}
\label{eq:thr}
  \begin{split}
   v \leq |q|/\sqrt 3 \ & \Longrightarrow \ \text{no transmitted soliton $u_T$,} \\
   v \geq \sqrt 3 |q|  \ & \Longrightarrow \ \text{no reflected soliton $u_R$.} 
\end{split} \end{equation}

Scattering of solitons by delta impurities is a natural model explored extensively in the physics literature -- see for instance \cite{CM},\cite{GHW}, and references given there.  The heuristic insight that at high velocities ``linear scattering'' by the external potential should dominate the partition of mass is certainly present there.  It would be interesting to see if bright solitons seen in Bose-Einstein condensates \cite{ASHSP} could be ``split'' using lasers modeled by delta impurities\footnote{That this might be related to the topic of this paper was suggested to the authors by N. Berloff.}.

In the mathematical literature the dynamics of solitons in the presence of external potentials has been studied in high velocity or semiclassical limits following the work of Floer and Weinstein \cite{FlWe}, and Bronski and Jerrard \cite{BJ}.  Roughly speaking, the soliton evolves according to the classical motion of a particle in the external potential. That is similar to the phenomena in other settings, such as the motion of the Landau-Ginzburg vortices.

The possible novelty in \eqref{eq:num1} and \eqref{eq:th2} lies in seeing {\em quantum} effects of the external potential strongly affecting soliton dynamics. 
The rest of the paper is organized as follows. In \S \ref{nr} we present 
a more detailed discussion of numerical results, and in \S \ref{nm} we 
outline the methods used in our computations.
Finally in \S \ref{ros} we discuss weaker but mathematically rigorous 
versions of \eqref{eq:num1} and \eqref{eq:th2} and give basic ideas 
behind the proofs.

\section{Numerical results}
\label{nr}
We now give numerical evidence for the results presented in 
\S \ref{in}: the asymptotics 
\eqref{eq:num1} and \eqref{eq:th2}.
We stress that the rigorous results of \cite{HMZ} provide
weaker error estimates and hold in a limited time range only.
We find it very interesting however that some results which are 
theoretically demonstrated, such as the thresholds \eqref{eq:thr}
or the long time behaviour of the the free NLS, are difficult to 
verify numerically -- see \S\S \ref{resow} and \ref{cnls}. 
On the other hand things which are hard to prove, such as the
existence of limits \eqref{eq:tqs}, seem to be very clear
in numerical analysis.

\begin{figure}
\scalebox{0.7}{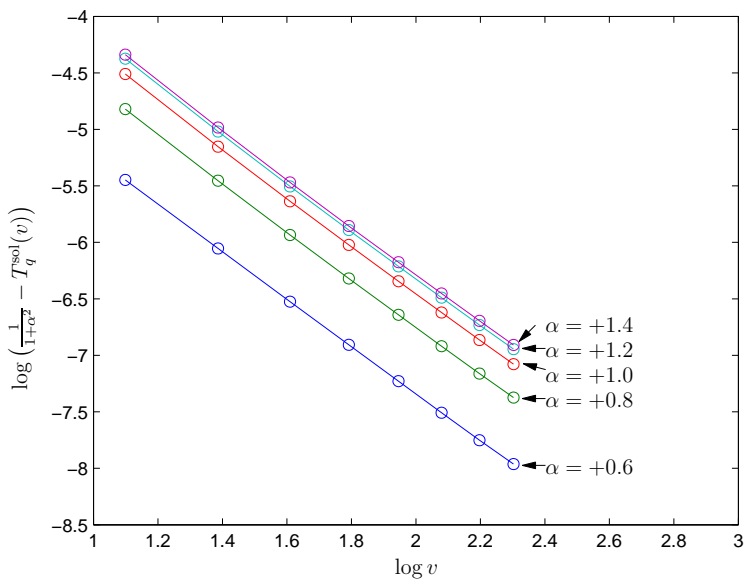}
\caption{This graph is a log--log plot showing that measured values of 
$T_q^{\rm{s}}(v)$ for $v=3,4,\ldots, 10$ conform well to the asymptotic formula \eqref{E:trans}. }
\label{F:trans}
\end{figure}

\begin{figure}
\scalebox{0.7}{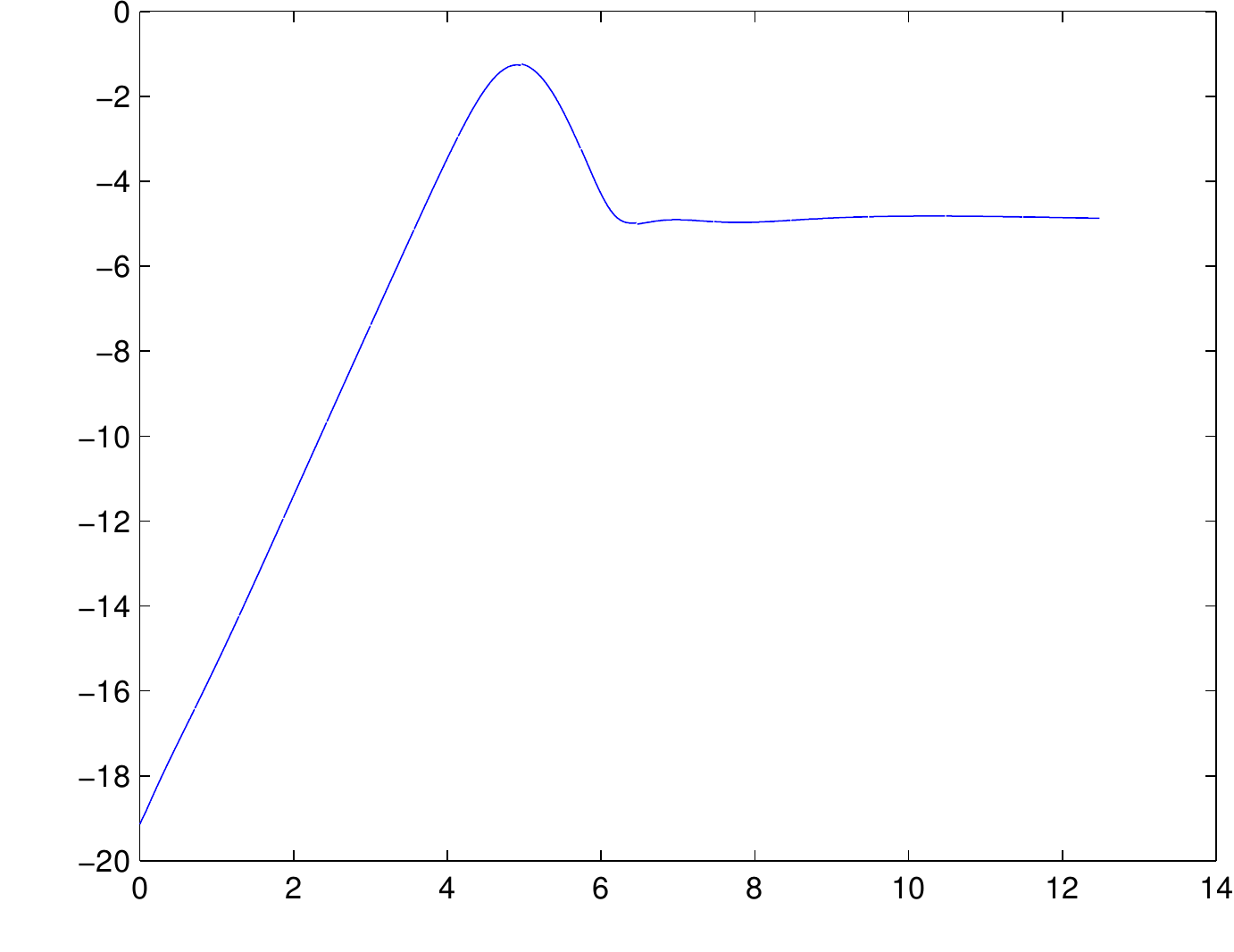}
\scalebox{0.7}{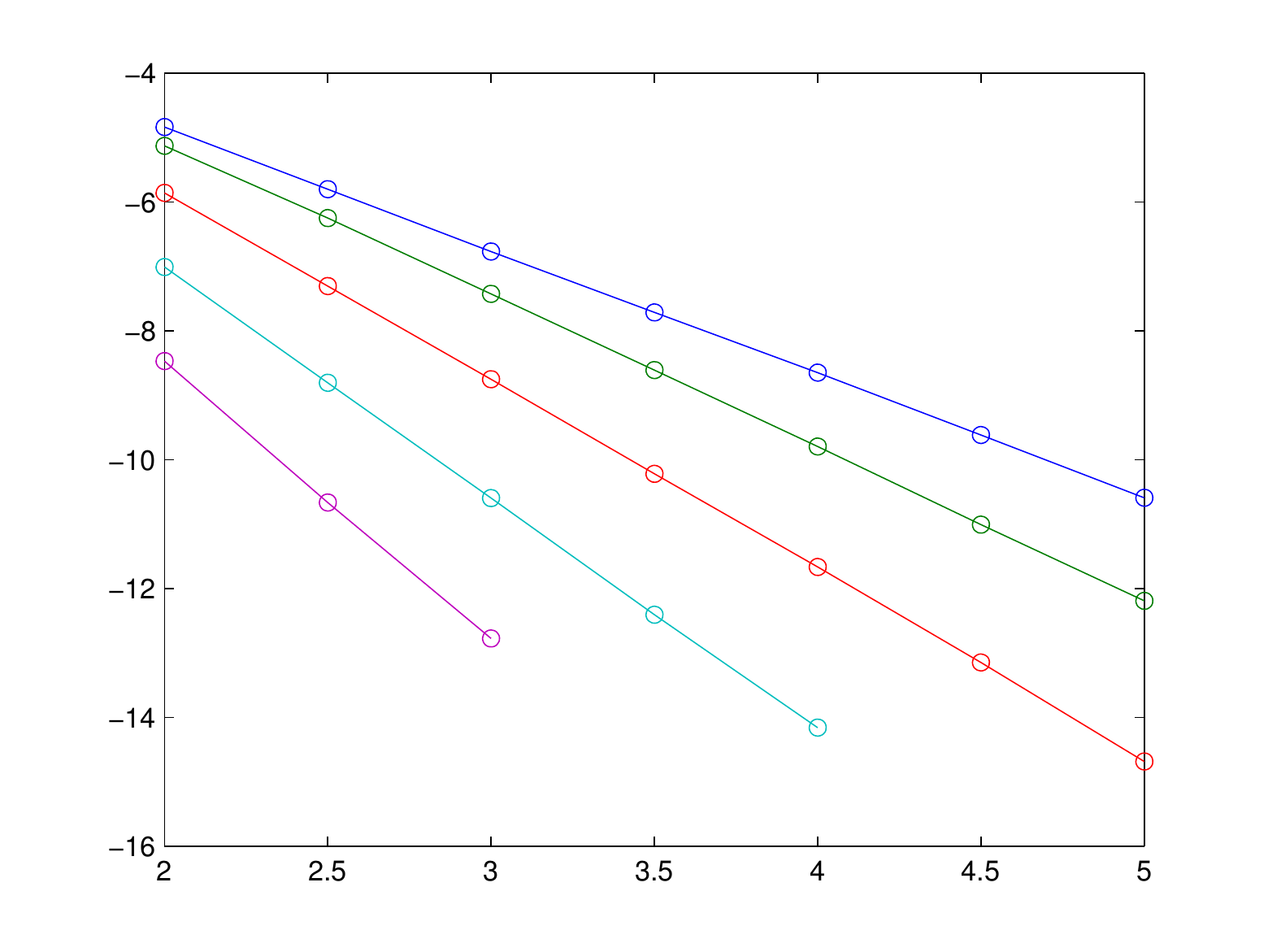}
\caption{The first plot shows, for $\alpha=-0.6$, $v=2$, $x_0=-10$, the stabilization of the value of $\int_{-0.5}^{0.5} |u(x,t)|^2 \, dx$ (which enters into the definition of $B_q^{\rm{s}}(v)$) after the interaction.  The second plot shows that the measured value of $B_q^{\rm{s}}(v)$ for velocities $v=2.0, 2.5, \ldots, 5.0$ conforms well to the asymptotic formula \eqref{E:trapped}. 
The values of $B_q^{\rm{s}}(v)$ that are $<e^{-14}$ cannot be measured with adequate precision, and thus the data for $|\alpha|\geq 1.2$ is limited.}
\label{F:trapped}
\end{figure}

\subsection{Asymptotics for the transmission and trapped rates}
We recall from \S \ref{in} the definitions of the transmission and reflection rates
\begin{align*}
T_q^{{{\rm{s}}}}(v) &= \frac{1}{2} \lim_{t \to +\infty} \int_{1/2}^{+\infty} |u(x,t)|^2 \,dx \\
R_q^{{{\rm{s}}}}(v) &= \frac{1}{2} \lim_{t \to +\infty} \int_{-\infty}^{-1/2} |u(x,t)|^2 \,dx \,.
\end{align*}
These definitions assume the existence of the limits. The numerical 
evidence strongly supports that the limits indeed exist.
We also define the trapped rate
\[ B_q^{{{\rm{s}}}}(v) = \frac{1}{2} \lim_{t \to +\infty} \int_{-1/2}^{+1/2} |u(x,t)|^2 \,dx \]
We have 
\[ B_q^{{{\rm{s}}}}(v)=0 \,, \ \ q \geq 0 \,. \]
However, $B_q^{{\rm{s}}}(v)>0$ for $q<0$ 
due to the presence of a bound state at $ \lambda = - i q $ for the 
linear operator $ H_q$:
\begin{equation}
\label{eq:eigen} \phi ( x ) = \sqrt{2 |q|} e^{ q |x| }  \,.
\end{equation}
The nonlinear problem has a bound state as well \cite{CM},\cite{GHW}:
\begin{equation}
\label{E:boundsol}
u(x,t) = e^{i\lambda^2t/2}\lambda \sech\left(\lambda|x| + \tanh^{-1}
\left({|q|}/{\lambda}\right) \right), \quad 0<\lambda < |q| \,.
\end{equation}
This bound state is ``left behind'' after the interaction. The parameter
$ \lambda $ depends on the initial condition.

\begin{table}
\begin{tabular}{|c|c|c|}
\hline
$\alpha$  & $a(\alpha)$ & $b(\alpha)$\\
\hline 
0.6 & 0.0415 & 2.0786 \\
0.8 & 0.0748 & 2.0788 \\
1.0 & 0.1007 & 2.0788 \\
1.2 & 0.1147 & 2.0778 \\
1.4 & 0.1185 & 2.0762 \\
\hline
\end{tabular} \qquad
\begin{tabular}{|c|c|c|}
\hline
$\alpha$  & $a(\alpha)$ & $b(\alpha)$\\
\hline 
-0.6 & 0.0441 & 2.1076 \\
-0.8 & 0.0761 & 2.0873 \\
-1.0 & 0.1014 & 2.0823 \\
-1.2 & 0.1151 & 2.0798 \\
-1.4 & 0.1189 & 2.0776 \\
\hline
\end{tabular}\qquad
\vspace{0.2in}
\caption{The left table gives the regression coefficients for \eqref{E:trans} for $q>0$, and the right table for \eqref{E:trans} for $q<0$ -- see 
\eqref{E:trans} for the definitions of $ a ( \alpha ) $ and $ b ( \alpha ) $,
$ \alpha = q/ v $.}
\label{T:regression}
\end{table}

In both the $q>0$ and $q<0$ cases, 
$T_q^{{\rm{s}}}(v)$ is numerically shown to follow an asymptotic
\begin{equation}
\label{E:trans}
T_q^{{{\rm{s}}}}(v) \sim \frac{1}{1+\alpha^2} - \frac{a(\alpha)}{v^{b(\alpha)}}
\,, \ \ \ \ \alpha = \frac{q}{v} \,. 
\end{equation}
Fig.\ref{F:trans} is a log--log plot of data  for velocities 3 to 10 
and $\alpha=+0.6$ to $+1.4$.  Results of linear regression on the subset of the data for velocities 5 to 10 yields the values of $a(\alpha)$ and $b(\alpha)$ reported in Table \ref{T:regression}.  We see that the value of $b(\alpha)$ displays little variation with $\alpha$ and is approximately $2.07$.  When the regression is performed on the subset of the data for velocities 8 to 10, the values of $b(\alpha)$ obtained are approximately $2.05$ leading us to conjecture that the true value of $b(\alpha)$ is exactly $2$, 
as stated in \eqref{eq:num1}.

  The data for $q<0$ gives a plot nearly identical to Fig.\ref{F:trans} owing to the exponential decay of the trapped mass $B_q^{{\rm{s}}}(v)$(discussed below).  The values of $a(\alpha)$ and $b(\alpha)$ obtained from the $q<0$ data are also reported in Table \ref{T:regression}, and we still expect that the true value of $b(\alpha)$ is $2$.  

Another feature apparent in Fig.\ref{F:trans} and Table \ref{T:regression} is that the value of $a(\alpha)$ stabilizes as $\alpha \to \infty$ (note the proximity of the lines for $\alpha=1.2$ and $\alpha=1.4$ in Fig.\ref{F:trans}).  This feature coincides with our analytical result, which establishes a (nonoptimal) bound on the asymptotic of $v^{-1/2+}$, independent of $\alpha$.   

\begin{table}
\begin{tabular}{|c|c|c|}
\hline
$\alpha$  & $d(\alpha)$ & $f(\alpha)$\\
\hline
-0.6 & 0.3610 & 1.9121 \\
-0.8 & 0.6952 & 2.3619 \\
-1.0 & 1.0328 & 2.9331 \\
-1.2 & 1.1506 & 3.5784 \\
-1.4 & 1.1351 & 4.3054 \\
\hline 
\end{tabular} 
\vspace{0.2in}
\caption{Numerical results for the parameters in 
\eqref{E:trapped} for $q<0$.}
\label{T:regression1}
\end{table}

The trapped mass coefficient $B_q^{{\rm{s}}}(v)$, on the other hand, decays exponentially. 
\begin{equation}
\label{E:trapped}
B_q^{{\rm{s}}}(v) \sim d(\alpha) e^{-f(\alpha)v} \,, \ \ 
\alpha = \frac{q}{v} \,. 
\end{equation}
The second frame of Fig.\ref{F:trapped}, which presents data for velocities 
2 to 5 and $\alpha=-0.6$ to $-1.4$, demonstrates that $B_q^{{\rm{s}}}(v)$ conforms well to the formula \eqref{E:trapped}.  Linear regression on the data yields the values of $d(\alpha)$ and $f(\alpha)$ reported in 
Table \ref{T:regression1}. 

 We see that in constrast to the behaviour of $b(\alpha)$, $f(\alpha)$ 
increases with $\alpha$.  
This produces a numerical road block 
in studying the asymptotic behaviour further 
for $\alpha=-1.2$ and $\alpha=-1.4$. We do not have 
enough significant digits in our data to measure 
values of $B_q^{{\rm{s}}}(v)$ less than $<e^{-14}$. 

For the nonlinear bound state, $u(x,t)$, given by  \eqref{E:boundsol}, we 
have 
$$\|u(\cdot, t)\|_{L^2}^2 = 2 (\lambda - {|q|} )\,. $$
Hence, the behaviour of $ B_q^{\rm{s}} ( v ) $ 
(see \eqref{E:trapped} and Table \ref{T:regression1}) shows 
that $\lambda$ approaches $q$ at an exponential rate as $v\to \infty$.
We note that for $\lambda$ very close to $|q|$, \eqref{E:boundsol} is approximately\footnote{This is obtained using the two approximations $\tanh x \approx 1-2e^{-2x}$ and $\sech x \approx 2e^{-x}$ for $x$ large.}
$$u(x,t) \approx e^{i\lambda^2t/2}\sqrt 2 \left( 1- \frac{|q|}{\lambda} \right) \lambda e^{-\lambda |x|} \approx e^{iq^2t/2} \sqrt 2 \left( 1- \frac{|q|}{\lambda} \right) |q| e^{- |qx|}$$
which is a multiple of the eigenstate \eqref{eq:eigen}.

Given that the interaction with the delta potential is dominated 
by the linear part of the equation, 
we expect the trapped state will, immediately after interaction, 
resemble a linear eigenstate that will then resolve on a longer time scale to a nonlinear bound state of the form \eqref{E:boundsol}.  
It is thus perhaps more appropriate to reverse the heuristics:  
Given a very small amplitude $A$, if we set 
\[ \lambda = |q|\left(1-\frac{A}{\sqrt 2}\right)^{-1} \,, \]
 we obtain that the eigenstate of $ H_q $ with amplitude $A$ is close to 
$u(x,t)$ given by \eqref{E:boundsol}, that is,
\begin{equation}
\label{E:boundapprox}
Ae^{iq^2t/2}e^{-|qx|} \approx u(x,t) \text{ solving \eqref{E:boundsol}}
\end{equation}
It is reasonable to expect that the nonlinear bound state ultimately selected from an immediate post-interaction eigenstate of the {\em linear} operator
$ H_q $, is ``close'', 
in the sense of \eqref{E:boundapprox}, 
to the starting eigenstate of $ H_q $, and indeed, the numerics point in this direction.  The first frame of Fig.\ref{F:trapped} shows for a typical case ($\alpha=-0.6$, $v=2$) that a stable trapped mass is selected within a reasonable amount of time following the interaction and there is little evidence of mass being radiated away from the origin.

\subsection{Resolution of outgoing waves}
\label{resow}
The stabilization of solitons described in \eqref{eq:th2} (and in 
a slightly weaker form rigorously in \cite[Theorem 2]{HMZ}) occurs 
over long time intervals -- see the comment at the end of \S \ref{ros}.

\begin{figure}
\scalebox{0.6}{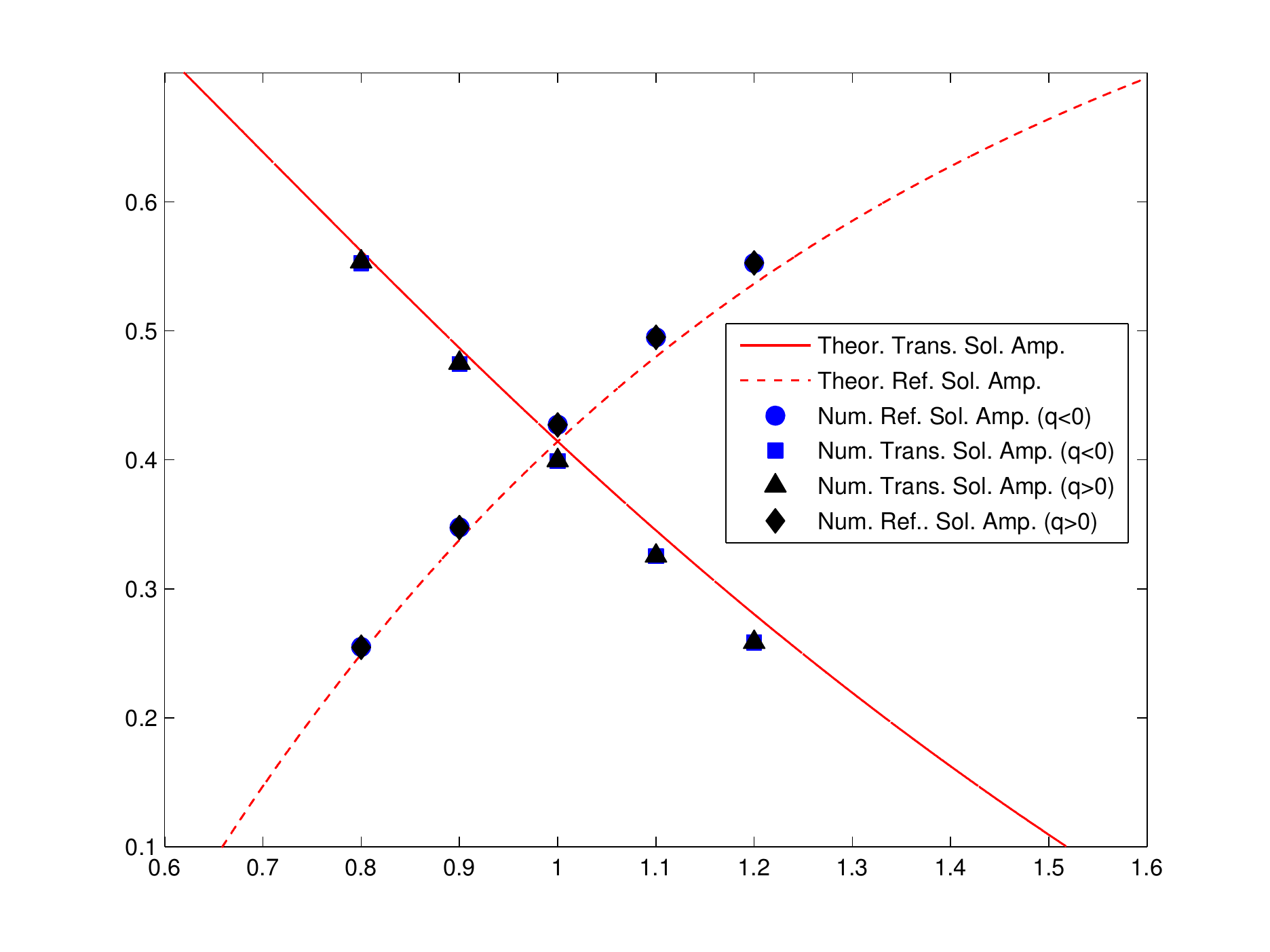}
\scalebox{0.6}{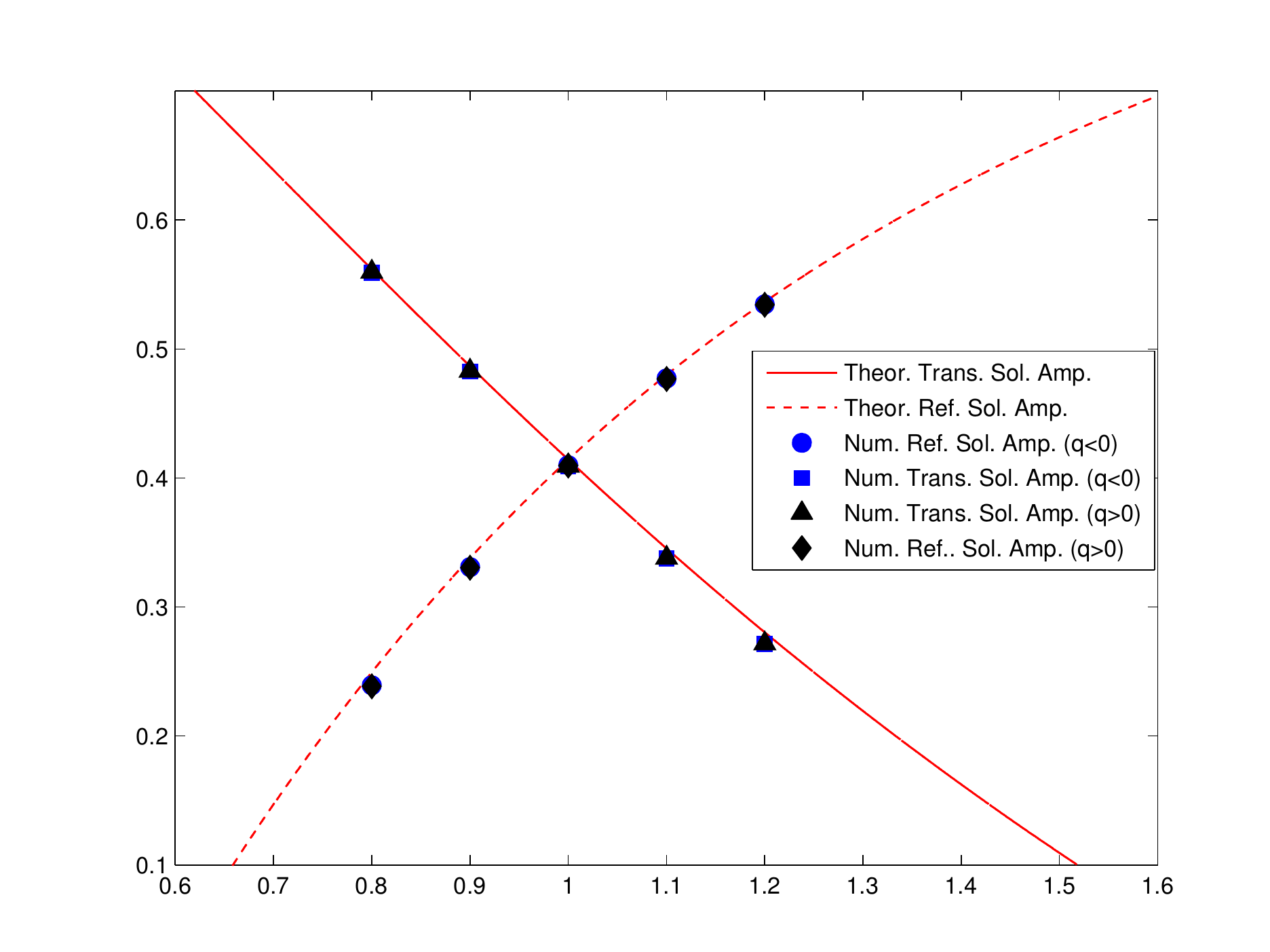}
\caption{A comparison of theory with numerically computed soliton parameters $ A_R $ and $ A_T $.  The first plot is for velocity 3; the second, which shows a better agreement, is for velocity 10.   Note here we are plotting \textit{amplitudes}, as opposed to fraction of the total \textit{mass} as in Fig.\ref{f:new}.}
\label{F:Maciejfav}
\end{figure}

Hence for the calculation of the amplitudes $ A_T $ and $ A_R $ 
in \eqref{eq:th2} we must alter our approach.  
We begin by solving the nonlinear equation \eqref{eq:nls}.  
However, we also measure the $L^2$ difference between the solution 
at every time step and the expected profile given by
\begin{eqnarray*}
\exp({ivx+i t(1-v^2)/{2}})  [t_q(v) \sech(x-x_0-vt) + r_q(v) \sech(x+x_0+vt)].
\end{eqnarray*}
Shortly after the time of interaction with the delta potential,
we see this difference attains a minimum. 
At this time, we save the computed solution, $ u ( x) $, and continue 
to solve forward in time.  As we solve forward, we compute
\[ \| \nlsq(t) u - \nlso(t)u \|_{L^2_x } \,.\]
Since this norm remains negligible we switch to the analysis of the
simpler solution of the \eqref{eq:nls1}, $ \nlso(t) u $.

More precisely, we truncate the solution on both sides of the delta to 
give two sets of initial data.  
We then perform a phase shift in order to give 
each piece a zero velocity. They are then embedded at the center of 
a very large grid with zeroes outside their computed ranges.  
From here, we solve forward on this larger grid using 
$\nlso(t)$ in order to observe the amplitude stabilization predicted by 
\eqref{eq:th2}. The grid is chosen large enough so that several 
amplitude oscillations can occur without interference 
from accumulated errors at the boundary.  
Though the amplitude continues to oscillate, we time average the 
amplitudes until we see stabilization.  

It is these time averages over significantly large intervals
that are reported
in  Fig.\ref{F:Maciejfav}.
We see a clear agreement with \eqref{eq:th2} especially for the higher
velocity. However, the theoretically predicted thresholds for the formation of
the reflected and transmitted solitons \eqref{eq:thr} are hard
to verify numerically.

\subsection{Confirmation of the free NLS asymptotics for initial data $\alpha \sech x$}
\label{cnls}

\begin{figure}
\scalebox{0.7}{\includegraphics{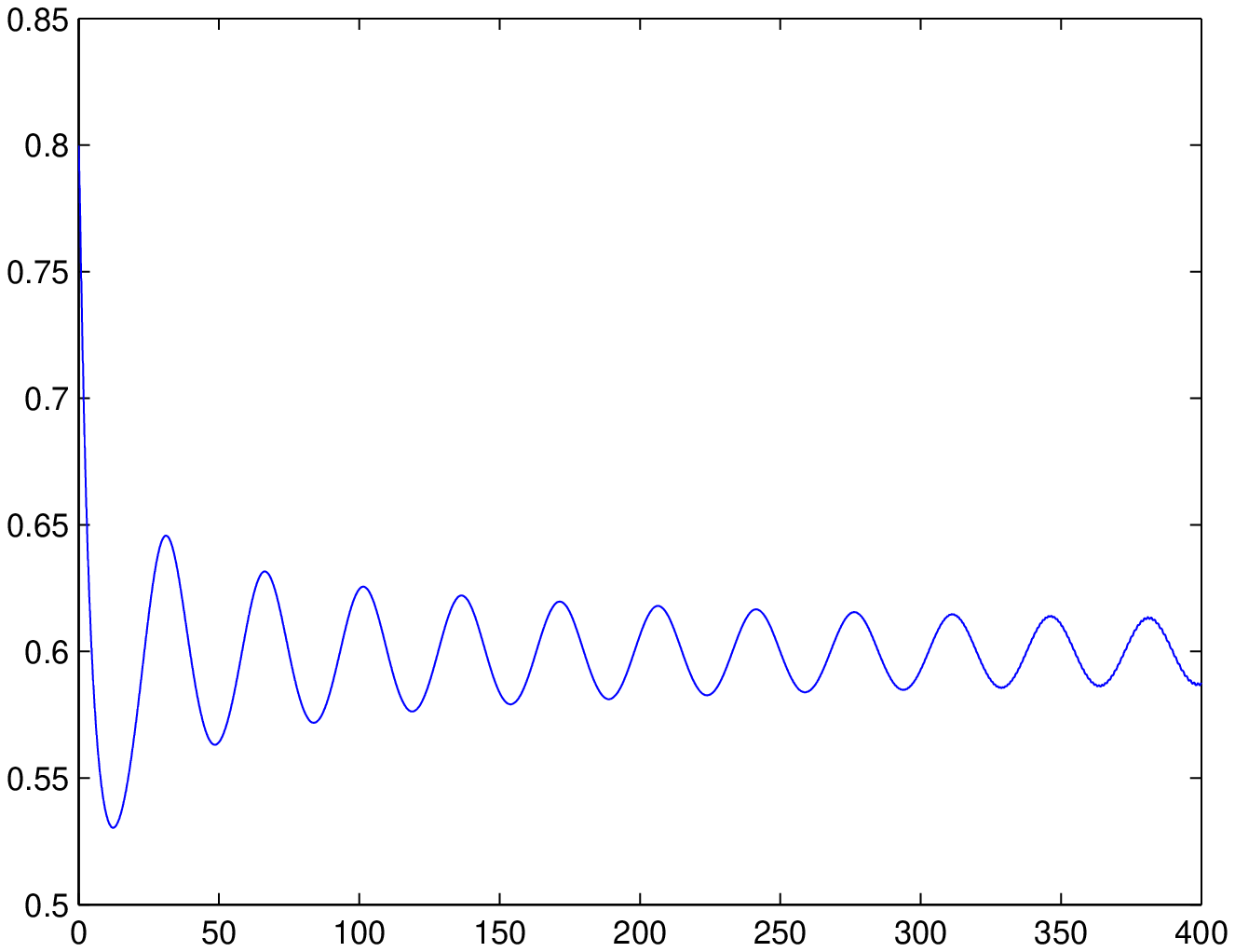}}
\scalebox{0.7}{\includegraphics{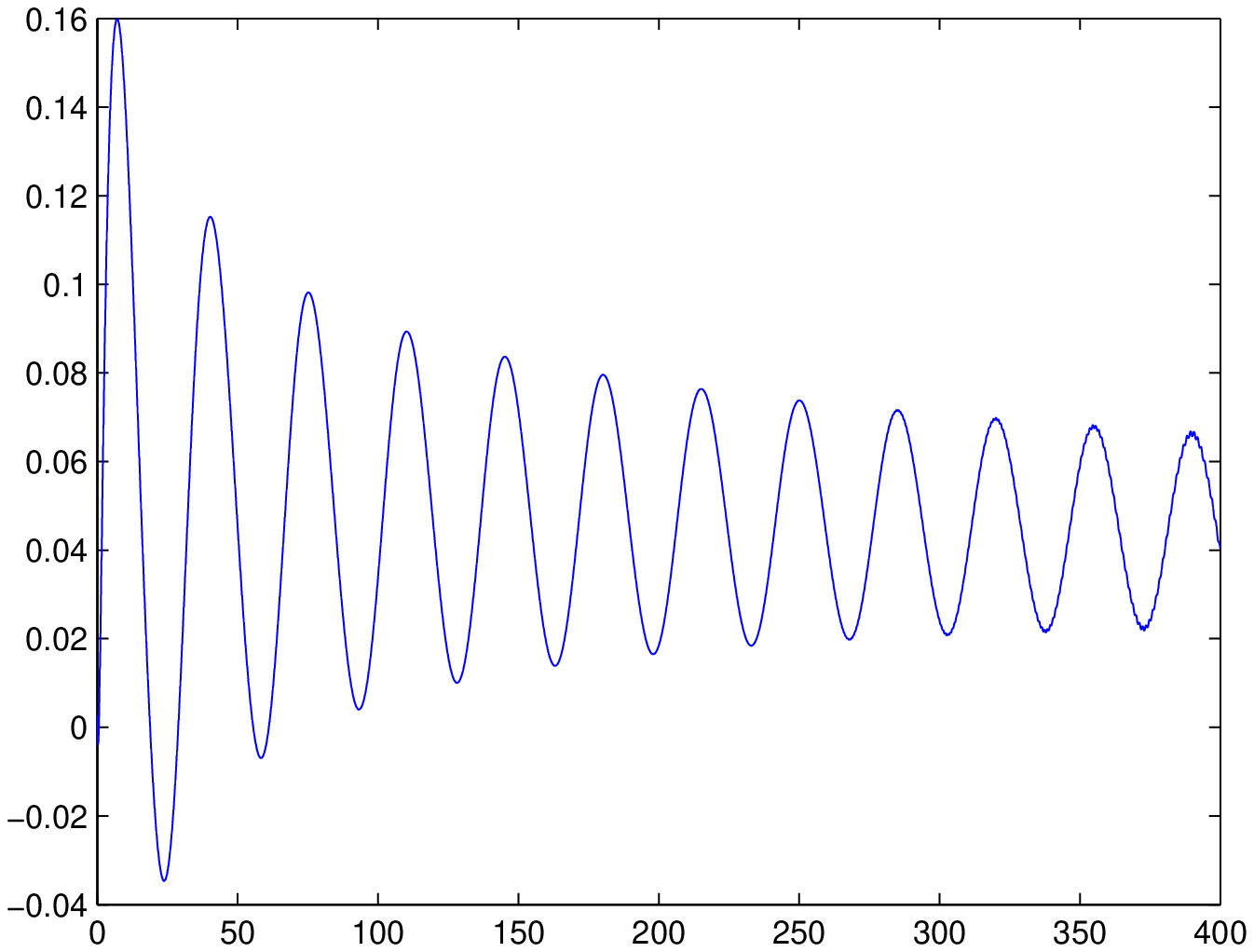}}
\caption{Two plots depicting the amplitude and phase of $u(0,t)$ for the free equation with initial data $0.8 \sech x$.}
\label{F:resolve}
\end{figure}

We now turn to the matter of propagating initial data $\alpha \sech x$ according to $i\partial_t u + \frac{1}{2}\partial_x^2 u + |u|^2u$, which has been explored analytically via the inverse scattering method in 
\cite[Appendix B]{HMZ}. 
Fig.\ref{F:resolve} reports the results of an experiment with 
$\alpha =0.8$.  The first panel depicts the time evolution of 
the amplitude at the spatial origin, $|u(0,t)|$, and 
the second panel depicts the deviation of the time evolution 
of the phase at the spatial origin from that of the soliton $(2\alpha-1)e^{i(2\alpha-1)^2t/2}\sech((2\alpha-1)x)$.  
The amplitude appears to be converging to the theoretically 
expected value of $2\alpha-1=0.6$ and the phase deviation appears to be converging to the expected value of $ \varphi ( 0.8) = 0.045$. 

In regard to the phase computation, it should be noted that although this experiment was performed on a $(x,t)$ grid of size $15000\times 20000$ with spatial extent $-600\leq x\leq 600$, the reported phase deviation is the difference of two numbers on the order of $100$ and the obtained values are three orders of magnitude smaller.  
This opens a possibility of an inaccuracy of this long-time computation. 
 
\section{Numerical methods}
\label{nm}
In this section we outline numerical methods used to produce the 
results described in \S \ref{nr}.
We discretize our equation,
\begin{eqnarray*}
&i u_t + u_{xx} + |u|^2 u - q \delta_0 (x) u = 0, \\
&u(0,x) = u_0,
\end{eqnarray*}
using a finite element scheme in space and the standard midpoint rule in time. Just as the equation itself this method is $L^2$ conservative.  
A finite difference scheme 
can also be implemented with an approximate delta function, but 
then convergence must be determined in terms of finer meshes as well 
as more accurate delta function approximations.  
In finite element methods the inherent integration allows us to directly incorporate the delta function into the discretization of the problem.  
A similar scheme was implemented without potential in {\cite{ADKM3}}, where 
the blow-up for NLS in several dimensions was analyzed.  

We review the method here and describe how to discretize the delta potential.  
Note that we require the spatial grid to be large enough to insure 
negligible interaction with the boundary.  
For the convergence of such methods without potentials see the
references in \cite{ADKM3}.


We select a symmetric region about the origin, $[-R,R]$, upon which 
we place a mesh of $N$ elements.  The standard hat function basis is 
used in the Galerkin approximation. We allow for a finer grid in a 
neighbourhood of length $1$ centered at the origin to better study the effects of the interaction with the delta potential.  
In terms of the hat basis the problem becomes:
\begin{eqnarray*}
&\langle u_t, v \rangle + {i} \langle u_x, v_x \rangle /2  
- i \langle |u|^2 u,v \rangle + i q u(0) v(0) = 0, \\
&u(0,x) = u_0 \,, \ \  u ( t , x) = \sum_{v} c_v ( t) v \,,
\end{eqnarray*} 
where $\langle \cdot, \cdot \rangle$ is the standard $L^2$ inner product, $v$ is a basis function and $u$, $u_0$ are linear combinations of the $v$'s. 
We remark that since $ v$'s are continuous the pairing of $  u v  $ with
the delta function is justified.

Since the $v$'s are hat functions, we have created 
a tridiagonal linear system with one contribution to the central element 
resulting from the delta function.  
Let $h_t > 0 $ be a uniform time step, and let 
\[ u_n = \sum_{v} c_v ( n h_t ) v \,, \]
be the approximate solution at the $n$th time step.
Implementing the midpoint rule in time, the system becomes:
\[ \begin{split}
& \langle u_{n+1} - u_n , v\rangle + i h_t \left\langle \left( {(u_{n+1} + u_{n})}/{2} \right)_x, v_x \right\rangle  + i h_t q ({u_{n+1}(0) + u_{n}(0)}) v(0)/2 \\
& \ \ \ \ = 
 i h_t \left\langle |({u_{n+1} + u_{n}})/{2}|^2 ({u_{n+1} + u_{n}})/{2} , v \right\rangle, \ \ \ \
 u_0 = \sum_v \alpha_v v,
\end{split} \]
By defining 
\[ y_{n} = (u_{n+1} + u_n)/2 \,, \] 
we have simplified our system to:
\[
\langle y_n , v \rangle + i \frac{h_t}{4} \langle ({y_{n}})_x,v_{x} \rangle  +  i \frac{h_t q}{2} y_n (0) v(0) = 
i \frac{h_t}{2} \langle |y_n|^2 y_n, v \rangle  +  \langle u_n, v\rangle.
\]
An iteration method from {\cite{ADKM3}} is now 
used to solve this nonlinear system of equations.  
To wit, 
\[ 
\langle y^{k+1}_{n} , v \rangle +
i \frac{h_t}{4} \langle (y^{k+1}_n)_x, v_{x} \rangle  + 
 i \frac{h_t q}{2} y^{k+1}_n (0) v(0) =
i \frac{h_t}{2} \langle |y_n^k|^2 y_n^k, v \rangle  +  \langle u_n, v \rangle.
\]
We take $y_n^0 = u_n$ 
and perform three iterations in order to obtain an approximate solution. 

The measure of success of this method is indicated by the agreement 
of the computed solutions and the exact solutions
such as \eqref{eq:trwv}, and more remarkably by 
the agreement with the inverse scattering method -- see \S \ref{cnls} above.
When the delta potential is present ($q\neq 0$) the $ L^2 $ norm 
remains essentially constant as predicted by theory.

\section{Review of theoretical methods}
\label{ros}

The theoretical results are proved in \cite{HMZ} in the 
case of a positive coupling constant, that is for $ q > 0 $. 
The asymptotic formul{\ae} \eqref{eq:num1} and \eqref{eq:th2}
are obtained rigorously for times
\begin{equation}
\label{eq:time}
\frac{|x_0|}{v} + v^{-1+\epsilon} \leq t \leq  \epsilon \log v\,, 
\end{equation}
and with the error term $ {\mathcal O}_{L^2} ( v^{-2} ) $ replaced by 
$ {\mathcal O}_{L^2} ( v^{ -( 1 -3 \epsilon)/2 }  ) $ for any $ \epsilon > 0 $,
provided that $ x_0 $, the center of the initial soliton, 
satisfies 
\[  x_0 < - v^{\epsilon } \,. \]

The starting point of the argument
is the observation
that the {\em Strichartz estimates} for $ H_q $ hold uniformly 
in $ q \geq 0 $. Strichartz estimates for dispersive equations \cite{SE}
describe joint space-time decay properties of solutions 
and are crucial in the control of interaction terms in nonlinear
equations. 

More precisely, for the problem
\begin{equation}
\label{eq:eq}
 i \partial_t u ( x , t ) + \tfrac{1}{2}\partial_{x}^2 u ( x, t ) 
- q \delta_0 ( x ) u ( x , t ) = f ( x , t ) \,, \ \ u ( x , 0 ) = u_0 ( x ) 
\,.\end{equation}
they generalize the well known energy inequality:
\[  \| u \|_{L^\infty_t L^2_x } \leq C \| u_0 \|_{L^2} + 
C \| f \|_{L^1_t L^2_x} \,,\]
where, for $ p \neq \infty $,
\[ \| u \|_{L^p_t L^r_x} = \left( \int \left( \int | u ( x , t) |^r dx 
\right)^{\frac{p}{r}} dt \right)^{\frac 1 p } \,.\] 
The Strichartz estimates allow more general exponents $ p $ and $ r $ 
and, which is essential to us, they hold
with constants $ C = C (p, r , \tilde p , \tilde r ) $ independent of $ q $:
\begin{gather}
\label{eq:Str}
\begin{gathered}
\| u \|_{ L^p_t L^r_x } \leq C \| u_0 \|_{L^2 } + C \| f \|_{ L_t^{\tilde 
p} L_x^{\tilde r} } \,, \\ 
2 \leq p, r \leq \infty \,, \ \ 1 \leq \tilde p , \tilde r \leq 2 \,, \ \
\frac 2 p + \frac 1 r = \frac 12 \,, \ \ \ \frac 2 {\tilde p }
+ \frac 1 {\tilde r} = \frac 52 \,, 
\end{gathered}
\end{gather}
see \cite[Proposition 2.2]{HMZ}. 

We also recall how the reflection and transmission coefficients
defined in \S \ref{in} for stationary scattering enter in 
time evolutions: for smooth 
$\psi$ vanishing outside of $ [-b, -a] $, $ 0 < a < b $, we have 
\begin{equation*}
\begin{split}
&  e^{-itH_q}[e^{ixv}\psi(x)](x) = \\
& \ \ \ 
\left\{ \begin{array}{ll} 
r(v)e^{-itH_0}[e^{-ixv}\psi(-x)](x) + e^{-itH_0}[e^{ixv}\psi(x)](x) +
e ( x, t )\,,  & x < 0 \,, \\
t(v)e^{-itH_0}[e^{ixv}\psi(x)](x) + e ( x , t )\,,  & x > 0 \,, 
\end{array} \right. 
\end{split} \end{equation*}
where 
$$\|e(x,t)\|_{L_x^2} \leq \frac{1}{v}\|\partial_x \psi\|_{L^2}$$
uniformly in $t$ -- see \cite[Lemma 2.4]{HMZ}.

To describe the proof of weaker versions of 
\eqref{eq:num1} and \eqref{eq:th2}
it will be useful to denote by 
\[ \nlsq(t)\varphi( x )  = u ( x , t) \]
the solution to 
\[i\partial_tu + \tfrac{1}{2}\partial_x^2 u - q\delta_0(x)u + |u|^2u=0 \,, \ \
u ( x , 0 ) = \varphi( x ) \,.\]
When $ q > 0 $ we refer to $ \nlsq(t) \varphi $ as 
the ``perturbed nonlinear flow'' and when $ q = 0 $, as 
the ``free nonlinear flow''. Similarly $ \exp ( - i t H_q ) \varphi $ for 
$ q > 0 $ is the ``perturbed linear flow'', and $ \exp ( - i t H_0 ) 
\varphi $ is the ``free linear flow''.

As discussed in Sect.\ref{in} we are interested in 
\[  u(x,t) = \nlsq(t)u_0(x) \,, \ \ 
u_0(x) = e^{ixv}\sech(x-x_0) \,, \ \ v \gg 1 \,, \ \  x_0 \leq -v^\epsilon 
\,, \ 0<\epsilon<1 \,.\]

The proof of \eqref{eq:num1} and \eqref{eq:th2} (in the time 
interval \eqref{eq:time} and with the error term $ {\mathcal O}_{L^2} ( v^{-2} ) $ replaced by $ {\mathcal O}_{L^2} ( v^{ -( 1 - 3 \epsilon)/2 } ) $) 
procceeds
in four phases.

\noindent\textbf{Phase 1 (Pre-interaction)}.  Consider $0\leq t \leq t_1$, where $t_1 = |x_0|/v - v^{-1+\epsilon}$ so that 
$x_0+vt_1=-v^{\epsilon}$.  
The soliton has not yet encountered the delta obstacle 
and propagates according to the free nonlinear flow
\begin{equation}
\label{E:approx1}
u(x,t) = e^{-itv^2/2}e^{it/2}e^{ixv}\sech(x-x_0-vt) 
+ \mathcal{O}(qe^{-v^{\epsilon}}), \quad 0\leq t\leq t_1\,. 
\end{equation}
The analysis here is valid provided $v$ is greater than some 
absolute threshold (independent of $q$, $v$, or $\epsilon$).  
But if we further require that $v$ be sufficiently large so that 
\[ v^{-3/2}e^{v^{\epsilon}} \geq q/v \,, \]
then 
\[qe^{-v^{\epsilon}} \leq v^{-1/2}\leq v^{-(1-\epsilon)/2} \,. \]
This is the error that arises in the main argument of Phase 2 below.

\noindent\textbf{Phase 2 (Interaction)}.  
Let $t_2 = t_1+2v^{-1+\epsilon}$ and consider $t_1\leq t \leq t_2$.  
The incident soliton, beginning at position $-v^{\epsilon}$, 
encounters the delta obstacle and splits into a transmitted component 
and a reflected component, which by time $t=t_2$, are concentrated at 
positions $v^{\epsilon}$ and $-v^{\epsilon}$, respectively.  
More precisely, at the conclusion of this phase (at $t=t_2$),
\begin{equation}
\label{E:approx4}
u(x,t_2) = 
\begin{aligned}[t]
&t(v)e^{-it_2v^2/2}e^{it_2/2}e^{ixv}\sech(x-x_0-vt_2)\\
&+r(v)e^{-it_2v^2/2}e^{it_2/2}e^{-ixv}\sech(x+x_0+vt_2) \\
&+ \mathcal{O}(v^{-\frac{1}{2}( 1 - \epsilon) })
\end{aligned}
\end{equation}

This is the most interesting phase of the argument, which proceeds 
by using the following three observations
\begin{itemize}
\item The perturbed nonlinear flow is approximated by the perturbed linear flow for $t_1\leq t \leq t_2$.
\item The perturbed linear flow is split as the sum of a transmitted component and a reflected component, each expressed in terms of the free linear flow of soliton-like waveforms.
\item The free linear flow is approximated by the free nonlinear flow on $t_1\leq t \leq t_2$.  Thus, the soliton-like form of the transmitted and reflected components obtained above is preserved.
\end{itemize}
The brevity of the time interval 
$[t_1,t_2]$ is critical to the argument, 
and validates the approximation of linear flows by nonlinear flows. It is 
here that we used the independence of $ q > 0 $ in \eqref{eq:Str}.

\noindent\textbf{Phase 3 (Post-interaction)}.  
Let $t_3=t_2+ \epsilon \log v$, and consider $[t_2,t_3]$.  
The transmitted and reflected waves essentially do not encounter 
the delta potential and propagate according to the free nonlinear flow, 
\begin{equation}
\label{E:post}
u(x,t) = 
\begin{aligned}[t]
& e^{-itv^2/2}e^{it_2/2}e^{ixv}\nlso(t-t_2)[t(v)\sech(x)](x-x_0-tv) \\
&+ e^{-itv^2/2}e^{it_2/2}e^{-ixv}\nlso(t-t_2)[r(v)\sech(x)](x+x_0+tv) \\
&+\mathcal{O}(v^{-\frac{1}{2}( 1 - 3\epsilon)}), \qquad t_2\leq t \leq t_3
\end{aligned}
\end{equation}
This is proved by a perturbative argument that enables us to evolve 
forward a time $ \epsilon \log v$ at the expense of enlarging the 
error by a multiplicative factor of $e^{ \epsilon \log v} = v^{\epsilon}$.  
The error thus goes from $v^{-( 1- \epsilon) /2}$ at $t=t_2$ to 
$v^{ - ( 1 - \epsilon)/2 + \epsilon }$ at $t=t_3$.

\noindent\textbf{Phase 4 (Soliton resolution)}.  
The last phase uses \eqref{E:post} and the following result based on 
inverse scattering method:
\begin{equation}
\label{eq:apb}
\nlso ( \alpha \, \sech ) =
\left\{ \begin{array}{ll} 
 e^{ i \varphi ( \alpha ) } \nlso ( ( 2 \alpha - 1 ) \sech ( ( 2 \alpha - 1)
\bullet ) ) + {\mathcal O}_{L^\infty } ( t^{-\frac12} ) 
& 1/2 < \alpha < 1 \,, \\
{\mathcal O}_{L^\infty } ( t^{-\frac12} ) &  0 < \alpha < 1/2 \,, 
\end{array} \right.
\end{equation}
where 
\[ \varphi ( \alpha ) 
= \int_0^\infty \log \left( 1 + \frac{ \sin^2 \pi \alpha }{ \cosh^2 
 \pi \zeta  } \right) \frac{ \zeta}{ \zeta^2 + ( 2 \alpha -1)^2 } 
d \zeta \,, \ \ 1/2 < \alpha < 1 \,, \]
see \cite[Appendix B]{HMZ} for the proof and references. 
The crucial part of the argument involves an evaluation 
of the transmission and reflection coefficients for the Zakharov-Shabat
system $ L \psi = \lambda \psi $, 
\[  L = -iJ\partial_x + iJ Q \,, \ \
Q = Q (t,  x ) = \begin{bmatrix} \ \ 0 & \alpha \, 
\sech x \\ - \alpha \, \sech x 
  &  \ 0 \end{bmatrix}\,, \ \
\qquad J = \begin{bmatrix} -1 & 0 \\ \ \ 0 & 1\end{bmatrix}\,. \]
That is done
by a well known computation \cite{M81},\cite[Sect.3.4]{Maib}
which reappears in many scattering theories, from the free
$S$-matrix
in automorphic scattering, to Eckart barriers in quantum chemistry.
We quote the results:
\[ t ( \lambda ) = \frac
{\Gamma ( \frac 12 + \alpha - i \lambda ) \Gamma ( \frac12 -
\alpha - i \lambda ) }{ \Gamma ( \frac12 - i \lambda ) } \,, \ \
 b ( \lambda ) = i \frac {\sin \pi \alpha }{ \cosh \pi \lambda } \,, \ \
r ( \lambda ) = { b ( \lambda ) } { t ( \lambda ) } \,. \]
The pole of $ t ( \lambda ) $ in $ \Im \lambda > 0 $ determines
the soliton appearing in \eqref{eq:apb}. We refer to 
\cite[Appendix B]{HMZ} for a detailed discussion. We only mention that
the long time asymptotics in the case of defocusing NLS \cite{DIZ} 
show that the error estimate $ {\mathcal O}_{L^\infty} ( t^{-1/2} ) $
is optimal and the stabilization to the soliton occurs over a long
time.

\medskip

\noindent
{\sc Acknowledgments.} We would like to thank V. Dougalis, 
J. Wilkening, J. Strain, 
and M. Weinstein for helpful discussions during the preparation of this paper.  The work of the first author was supported in part by an NSF postdoctoral fellowship, and that of the second and third author by NSF grants  DMS-0354539 and DMS-0200732.


\begin{thebibliography}{XX}

\bibitem{ASHSP}
U.~Al Khawaja, H.T.C.~Stoof, R.G.~Hulet, K.E.~Strecker, and 
G.B.~Partridge, 
{\em Bright Soliton Trains of Trapped Bose-Einstein Condensates,
}  Phys. Rev. Lett. {\bf 89}(2002), 200404-



\bibitem{ADKM3} G.D. Akrivis, V. A. Dougalis, O. A. Karakashian, and W. R. McKinney, {\em Numerical approximation of blow-up of radially symmetric solutions of the nonlinear Schr\"odinger equation,} SIAM J. Sci. Comput. {\bf 25}(2003), no. 1, 186-212.

\bibitem{Br} D. Braess, {\em Finite Elements: Theory, fast solvers, and applications in solid mechanics,} Cambridge University Press, 1997.



\bibitem{BJ} J. C. Bronski and R. L. Jerrard, 
{\em Soliton dynamics in a potential,} Math. Res. Lett. 
{\bf 7}(2000), 329-342. 

\bibitem{CM} X.D. Cao and B.A. Malomed, {\em Soliton-defect collisions in the nonlinear Schr\"odinger equation,}
Physics Letters A {\bf 206}(1995), 177--182.


\bibitem{SE} T.~Cazenave, {\em Semilinear Schrödinger equations.} 
Courant Lecture Notes in Mathematics, {\bf 10},
American Mathematical Society, Providence, 2003.


\bibitem{DIZ} P.A. Deift, A.R. Its, and X. Zhou,
{\em Long-time asymptotics for integrable nonlinear wave
equations}, in {\em Important developments in soliton theory}, 
181--204, Springer Ser. Nonlinear Dynam., Springer, Berlin, 1993.



\bibitem{FlWe} A. Floer and A. Weinstein, {\em Nonspreading wave packets for the cubic Schr\"odinger equation with a bounded potential}, 
J. Funct. Anal. {\bf 69}(1986), 397--408. 


\bibitem{GHW} R.H. Goodman, P.J. Holmes, and M.I. Weinstein, 
{\em Strong NLS soliton-defect interactions,}
Physica D {\bf 192}(2004), 215--248.

\bibitem{HMZ} J. Holmer, J. Marzuola, and M. Zworski, 
{\em Fast soliton scattering by delta impurities,}
{\tt math.AP/0602187}, preprint 2006.





\bibitem{KiMa} Y.S. Kivshar and B.A. Malomed,
{\em Dynamics of solitons in nearly integrable systems,} 
 Rev. Mod. Phys. {\bf 61}(1989), 763--915.

\bibitem{LL} L.D.~Landau and E.M.~Lifshitz, 
{\em Quantum Mechanics,} 3rd Edition.


\bibitem{Maib} A.I. Maimistov and M. Basharov,
{\em Nonlinear Optical Waves,} 
Fundamental Theories of Physics, {\bf 104}
Kluwer Academic Publishers, Dordrecht, Boston, London, 1999


\bibitem{M81} J.W.\ Miles, {\em An envelope soliton problem,} SIAM J. Appl. Math. {\bf 41} (1981), no. 2, 227--230.




\bibitem{ZS72} V.E.\ Zakharov and A.B.\ Shabat, {\em Exact theory of two-dimensional self-focusing and one-dimensional self-modulation of waves in nonlinear media,} Soviet Physics JETP {\bf 34} (1972), no. 1, 62--69.

\end{thebibliography}
\end{document}